\documentclass[12pt]{amsart}
\usepackage{latexsym}           %
\usepackage{amssymb}
\usepackage{epsfig}             %
\usepackage{amsfonts}

\newcommand{\IF}{\mathbb{F}}

\newcommand{\IP}{\mathbb{P}}                                     
                           
\newcommand{\IC}{\mathbb{C}}
\newcommand{\IZ}{\mathbb{Z}}

\newcommand{\M}{\mathcal{M}}

\newcommand{\cC}{\mathcal{C}}

\newcommand{\F}{\mathcal{F}}
\newcommand{\cN}{\mathcal{N}}
\newcommand{\cO}{O}            %

\newcommand{\g}{       \mathfrak{g}     }

\newcommand{\pf}{\begin{bpf}}

\newcommand{\pfms}{\begin{bpfms}}
\newcommand{\epf}{\end{bpf}\hfill$\square$\\}           %
\newcommand{\epfms}{\end{bpfms}\hfill$\square$\\}               %

\newcommand{\al}{\alpha}
\newcommand{\be}{\beta}
\newcommand{\ga}{\gamma}

\newcommand{\tr}{\text{\rm Tr}}
\newcommand{\Hom}{\text{\rm Hom}}

\newcommand{\End}{\text{\rm End}}
\newcommand{\diag}{{\text{\rm diag}}}

\newcommand{\spq}{/\!\!/}

\def\mapright#1{\smash{
        \mathop{\longrightarrow}\limits^{#1}}}

\theoremstyle{plain}
\newtheorem {hypo}{\bf\hspace{-\parindent}Hypothesis}
\newtheorem*{plainthm}{Theorem}

\newtheorem {lem}[hypo]{Lemma}%

\theoremstyle{definition}

\theoremstyle{remark}
\newtheorem {rmk}[hypo]{Remark}%

\begin{document}

\title{Painlev\'e Equations and Complex Reflections}
\author{Philip Boalch}

\begin{abstract}
We will explain how some new algebraic solutions of the sixth Painlev\'e
equation arise from complex reflection groups, thereby extending some results
of Hitchin and Dubrovin--Mazzocco for real reflection groups.
The problem of finding explicit formulae for these solutions will be 
addressed elsewhere.

\ 

\noindent
R{\tiny \'ESUM\'E}.\quad 
Nous expliquerons comment de nouvelles solutions alg\'ebriques de la 
sixi\`eme \'equation de Painlev\'e proviennent des groupes complexes de 
r\'eflexion,
prolongeant les r\'esultats de Hitchin et de Dubrovin--Mazzocco 
pour les groupes r\'eels de r\'eflexion. Le probl\`eme de trouver des formules 
explicites pour ces solutions sera trait\'e ailleurs. 
\end{abstract}

\maketitle

\renewcommand{\baselinestretch}{1.0}            %
\normalsize

\begin{section}{Introduction}

One of the themes of this article is
the study of monodromy actions in nonabelian cohomology.
In \cite{KPS01} Katzarkov--Pantev--Simpson searched for
fundamental group representations having
dense orbits under the monodromy action, whereas we will consider the opposite
extreme of finite orbits, in the simplest non-trivial case.
Namely we will look for finite orbits of the monodromy action on the set of
$SL_2(\IC)$ representations of the fundamental group of the
four-punctured sphere, or equivalently for
finite branching solutions of the sixth Painlev\'e equation.
Somewhat surprisingly these will arise from certain subgroups of $GL_3(\IC)$.

Another theme is that of extending results relating to
real reflection groups (or Weyl groups) to complex reflection groups. 
Various instances of this trend have appeared recently in the literature 
(cf. e.g. Brou\'e--Malle--Michel \cite{BMM99} or Totaro \cite{tot}), although
the instance here has quite a different flavour.

The sixth Painlev\'e equation (PVI):
$$\eta''=\left(\frac{1}{\eta}+\frac{1}{\eta-1}+\frac{1}{\eta-t}\right)\frac{(\eta')^2}{2}
-\left(\frac{1}{t}+\frac{1}{t-1}+\frac{1}{\eta-t}\right)\eta'  $$
$$\quad\qquad\qquad+\frac{\eta(\eta-1)(\eta-t)}{t^2(t-1)^2}\left(
\al+\frac{\be t}{\eta^2} + \frac{\gamma(t-1)}{(\eta-1)^2}+
\frac{\delta t(t-1)}{(\eta-t)^2}\right) $$
(where $\al,\be,\ga,\delta\in\IC$  are parameters)
for a (local) function $\eta(t): \IP^1\setminus\{0,1,\infty\}\to \IC $
has the Painlev\'e property, and so any such local
solution extends to a meromorphic function on the universal cover of the 
three-punctured sphere.
Recently Dubrovin--Mazzocco \cite{DubMaz00} 
have classified all the algebraic solutions of the one-parameter family 
$\text{PVI}_\mu$ of Painlev\'e six equations having parameters of the form
$(\al,\be,\ga,\delta)=((2\mu-1)^2,0,0,1)/2$ for 
$\mu\in\IC$ with $2\mu\not\in\IZ$. 
They found that such algebraic solutions
are precisely the finite branching solutions\footnote{
Clearly the algebraic solutions are finite branching, but the converse is not
obvious; one knows a finite branching solution lifts to a meromorphic function
on a punctured algebraic curve (finite over $\IP^1\setminus \{0,1,\infty\}$) 
and
must show this function is algebraic. } 
and
that, up to equivalence, they are in one-to-one correspondence with braid group
orbits of generating triples of reflection groups in three dimensional
Euclidean space; for the tetrahedral and octahedral groups
there is just one such braid
group orbit and for the icosahedral group there are three, so they obtain five
solutions altogether.
Some of these algebraic solutions have also been independently discovered by
Hitchin \cite{Hit-Poncelet, Hit-Ed}, and Hitchin's approach may be summarised
as follows:

It is well-known (cf. e.g. \cite{JM81}) that PVI is equivalent to the equations
for isomonodromic deformations of a linear system of Fuchsian 
differential equations
of the form
\begin{equation} \label{eqn: 2x2 linear}
\frac{d\Phi}{dz}= A(z)\Phi; \qquad A(z)=\sum_{i=1}^3\frac{A_i}{z-a_i}
\end{equation}
(where the $A_i$'s are $2\times 2$ matrices) as the collection of 
pole positions $(a_1,a_2,a_3)$ is varied in 
$\IC^3\setminus{\text{diagonals}}$.
Such isomonodromic deformations are governed by Schlesinger's equations:
\begin{equation} \label{eqn: 2x2 schles}
\frac{\partial A_i}{\partial a_j}=\frac{[A_i,A_j]}{a_i-a_j}\qquad 
\text{if } i\ne j, \text{ and}\qquad
\frac{\partial A_i}{\partial a_i}=-\sum_{j\ne
i}\frac{[A_i,A_j]}{a_i-a_j}
\end{equation}
and it is shown in \cite{JM81} how 
these are equivalent to PVI, where the time $t$ is
the cross-ratio of the four pole positions $(a_1,a_2,a_3,\infty)$ of 
\eqref{eqn: 2x2 linear}, and the four parameters $\al,\be,\ga,\delta$ are
essentially the differences of the eigenvalues of the four residues 
$A_1,A_2,A_3,A_4:=-A_1-A_2-A_3$
of $A(z)dz$. 
By observing that Schlesinger's equations preserve the adjoint orbit $O_i$
containing each $A_i$ and 
are invariant under overall conjugation of $(A_1,A_2,A_3,A_4)$ 
one sees that Schlesinger's equations amount to a flat connection on the
trivial fibre bundle 
\begin{equation}\label{eqn: rk 2+ bundle}
\M^*:=(O_1\times O_2\times O_3\times O_4)\spq G \times B \mapright{} B
\end{equation}
over $B:=\IC^3\setminus{\text{diagonals}}$, where the fibre 
$(O_1\times\cdots \times O_4)\spq G$ is the quotient of 
$$\{(A_1,A_2,A_3,A_4)\in O_1\times O_2\times O_3\times O_4 \bigl\vert \sum A_i
=0\}$$ by overall conjugation by $G=GL_2(\IC)$. 
(Generically this fibre is two dimensional, relating to the fact that PVI is
second order, and has a natural complex symplectic structure.) 
Now for each point $(a_1,a_2,a_3)$ of the base $B$ one can also 
consider the set 
$$\Hom_\cC(\pi_1(\IC\setminus\{a_i\}),G)/G$$ of conjugacy classes of 
representations of the 
fundamental group of the four-punctured sphere, where representations are
restricted to take the simple loop around $a_i$
into the conjugacy class $\cC_i:=\exp(2\pi\sqrt{-1} \cO_i)\subset G$ 
($i=1,\ldots,4, a_4=\infty$).
These spaces  of representations are also generically two dimensional and fit
together into a fibre bundle
$$M\mapright{} B.$$
Moreover this bundle $M$ has a {\em complete} flat connection (the
``isomonodromy connection'')  defined locally by identifying representations
taking the same values on a fixed  set of fundamental group generators.
The Schlesinger equations are the pullback of the isomonodromy connection
along the natural bundle map
$$\nu: \M^*\mapright{}M$$ 
defined by 
taking the systems \eqref{eqn: 2x2 linear} to their monodromy representations
(cf. \cite{Hit95long, smid}).

Thus one approach to finding
algebraic solutions to PVI is to start by finding 
finite branching solutions to Schlesinger's equation. In turn these correspond
to finite branching sections of the isomonodromy connection. 
However the isomonodromy connection is complete and so its branching amounts
to an action of the fundamental group of the base $\pi_1(B)$ (the pure 
three-string braid group) on a fibre
$\Hom_\cC(\pi_1(\IP^1\setminus\{a_i\}),G)/G$. 
This action extends to an action of the
full braid group $B_3$ on $\Hom(\pi_1(\IP^1\setminus\{a_i\}),G)$ which 
in turn comes from  the standard action of $B_3$ on triples of generators of
the free group $\IF_3\cong \pi_1(\IP^1\setminus\{a_i\})$ on three generators.
In particular if one chooses a representation $\rho\in 
\Hom(\pi_1(\IP^1\setminus\{a_i\}),G)$ whose image
$\rho(\pi_1(\IP^1\setminus\{a_i\}))$ in $G$ is a {\em finite} group, then one
knows immediately that the braid group orbit containing $\rho$ is finite.
(Clearly choosing such a representation is equivalent to choosing a triple of
elements of $G$ which generate a finite subgroup.)
The algebraic solutions of Hitchin were found in this way, 
starting with 
finite subgroups of $SU_2\subset G=GL_2(\IC)$.
Dubrovin--Mazzocco use a different procedure to obtain finite orbits of the
braid group but note (\cite{DubMaz00} Remark 0.2) 
that Hitchin's solutions are related to
theirs by a symmetry of PVI. 

In this paper we will use a different representation of PVI as an isomonodromy
equation, this time for certain rank {\em three} Fuchsian linear systems, 
and thereby obtain finite branching solutions of PVI from certain finite
subgroups of $GL_3(\IC)$. The main result is the following:

\begin{plainthm}
For each triple of generators of a three-dimensional 
complex
reflection group there is a finite branching solution of the sixth Painlev\'e
equation. 
\end{plainthm}

Thus we have found that the word ``complex'' may be added to the statement of
Dubrovin--Mazzocco \cite{DubMaz00} and so more solutions are obtained. 

We recall that Shephard--Todd \cite{Shep-Todd} 
have classified the finite groups generated by complex reflections and showed
that in three-dimensions, apart from the real reflection groups, there are four
irreducible complex reflection groups generated by triples of reflections,
of orders 336, 648, 1296 and 2160
respectively, as well as two infinite families $G(m,p,3), m\ge 3, p=1,m$ of
groups of orders $6m^3/p$. 
For $m=2$ and $p=1,2$ 
these would be the symmetry groups
of the octahedron and tetrahedron respectively.
(In general, 
for other $p$ dividing $m$, $G(m,p,3)$ is not generated by a triple of
reflections and so does not satisfy the hypotheses of the theorem.)

Note that we will not address here the further problem 
of finding explicit formulae for these solutions.
(Also we will not prove the
algebraicity of these solutions here but remark that, in the example 
we consider in Section 4, 
this may be proved directly  using Jimbo's work \cite{Jimbo82} 
on the asymptotics of PVI.
In general some modification of \cite{Jimbo82}
would be necessary to prove the
algebraicity of all the other solutions coming from complex reflection groups.)

The organisation of this paper is as follows.
In section 2 we will describe the three-dimensional Fuchsian systems we are
interested in and explain how to relate the Schlesinger equations for their 
isomonodromic deformations to the full four-parameter family of 
PVI equations.
Section 3 will then describe the braid group action on the corresponding space
of fundamental group representations.
Finally Section 4 describes an example of 
the finite branching solutions of PVI 
that arise from triples of generators of 
three-dimensional complex reflection groups. 

\ 

{\small
{\bf Acknowledgments.}\ \ 
This article is a simplified version\footnote{In particular the relation to
  Stokes multipliers and one of the
  original motivations (to better understand the braid group actions of
  \cite{DKP, bafi} for $GL_n$) are no longer apparent, but will be
  elucidated elsewhere.} of the authors talk at the conference 
in honour of Fr\'ed\'eric Pham, Nice 1-5 July 2002; the author
 is grateful to the
organisers for the invitation. 
Some inspiration for writing this up was provided by Y. Ohyama's survey talk
at RIMS Kyoto September 2001.
Part of this work was carried out at I.R.M.A Strasbourg, 
supported by the E.D.G.E Research Training Network
HPRN-CT-2000-00101.
}
\end{section}
\begin{section}{The rank three systems}

Let $V=\IC^3$, define 
 $G=GL(V)$ now and let $\g=\End(V)$ denote its Lie algebra.
Choose three non-integral complex numbers $\lambda_i\in \IC\setminus\IZ$ and
let $\cO_i\subset\g$ be the adjoint orbit of rank-one matrices having
trace $\lambda_i$, for $i=1,2,3$.
Thus $\cO_i$ is four dimensional and 
consists of matrices conjugate to $\diag(\lambda_i,0,0)$ 
and any element $B_i\in\cO_i$ maybe written in the form
$$B_i=e_i\otimes\al_i\qquad\text{where}\quad e_i\in V, \al_i\in V^* 
\text{ and }
\al_i(e_i)=\lambda_i.$$
Also choose a generic adjoint orbit $\cO_4\subset \g$ (which has dimension
six) and consider the space
$$\F:=(\cO_1\times \cO_2\times \cO_3\times \cO_4)\spq G $$
which is defined as the quotient of 
\begin{equation}\label{eqn: prequot}
\{(B_1,B_2,B_3,B_4)\in \cO_1\times\cdots\times \cO_4 \bigl\vert \sum B_i
=0\}
\end{equation}
 by overall conjugation by $G$.
Observe that $\F$ is of dimension two %
and so heuristically we would expect the equations for 
isomonodromic deformations of the Fuchsian systems  
\begin{equation} \label{eqn: 3x3 linear}
\frac{d\Phi}{dz}= B(z)\Phi; \qquad B(z)=\sum_{i=1}^3\frac{B_i}{z-a_i}
\end{equation}
(with $B_i\in \cO_i$)
to be equivalent to a Painlev\'e equation (i.e. to a second order equation
with the Painlev\'e property); this is indeed the case and we will now show we 
in fact get PVI.

Without loss of generality let us restrict the pole positions of 
\eqref{eqn: 3x3 linear} to be $(a_1,a_2,a_3)= (0,1,t)$ with
$t\in\IP^1\setminus\{0,1,\infty\}$. Then the Schlesinger equations for
isomonodromic deformations of \eqref{eqn: 3x3 linear} take the form
\begin{equation} \label{eqn: 3x3 schles}
\frac{dB_1}{dt}= \frac{[B_3,B_1]}{t},\quad
\frac{dB_2}{dt}= \frac{[B_3,B_2]}{t-1},\quad
\frac{dB_3}{dt}= \frac{[B_1,B_3]}{t}+\frac{[B_2,B_3]}{t-1}
\end{equation}
where the third  equation follows from the first two and the fact that
$B_1+B_2+B_3=-B_4$ is constant.
As before these equation descend to define a connection on the trivial bundle
$$\cN^*:= \F\times \IP^1\setminus\{0,1,\infty\}\to \IP^1\setminus\{0,1,\infty\}$$
with fibre $\F$.
Let us choose some coordinates on $\F$ and then rewrite Schlesinger's 
equations in terms of them.
\begin{lem}
The functions $x:=\tr(B_1B_3)$ and $y:=\tr(B_2B_3)$ are local coordinates
near a generic point of $\F$.
\end{lem}
\pf
Let $e_1,e_2,e_3$ be a basis of $V$.
If $\cO_4$ is generic then each $G$ orbit in \eqref{eqn: prequot}
contains a point having $B_i=e_i\otimes\al_i$ for some $\al_i\in V^*$, for
$i=1,2,3$, and any other such point of the $G$-orbit is of the form 
$(uB_1u^{-1},uB_2u^{-1},uB_3u^{-1},uB_4u^{-1})$ for some diagonal matrix $u$.
Thus if $b_{ij}=\al_i(e_j) = (-B_4)_{ij}$ then we see
$$\F\cong \{B_4\in\cO_4 \bigl\vert (B_4)_{ii}=-\lambda_i\text{ for $i=1,2,3$ }\}/T$$
where the diagonal torus $T\subset G$ acts by conjugation. 
Thus we just have to examine the action of $T$ on the off-diagonal entries of 
$B_4$.
Now, the $T$-invariant functions of the off-diagonal entries of $B_4$ 
are generated by the five functions
$$w=b_{12}b_{21},\  x=b_{13}b_{31},\  y=b_{23}b_{32}, \ 
p=b_{12}b_{23}b_{31}, \ q=b_{13}b_{32}b_{21},$$
and they satisfy the relation $wxy=pq$,
so that $\F$ is embedded in the subvariety of $\IC^5$ cut out by
this equation.
The further equations determining $\F$ (corresponding to fixing the eigenvalues
$B_4$) maybe be written in terms of $w,x,y,p,q$ as follows.
The sum of the eigenvalues is fixed since the diagonal part of $B_4$ is fixed.
The sum of the squares of the eigenvalues is
$$\tr (B_4^2) = \tr(B_1+B_2+B_3)^2 = \lambda_1^2 +\lambda_2^2 +\lambda_3^2 
+ 2(w+x+y)$$
and so fixing this amounts to a constraint of the form
$w=c-x-y$
for a constant $c$.
Similarly fixing the sum of the cubes of the eigenvalues leads to a constraint
of the form
$q=-p+ax+by+k$
for constants $a,b,k$.
Thus eliminating $w,q$ we see $\F$ is locally identified with the variety
in $\IC^3\ni(x,y,p)$ with equation
$$xy(x+y-c)=p(p -ax-by-k).$$
Thus fixing some generic values of $(x,y)$
determines $p$ up to a sign, and so
$(x,y)$ are indeed generically good coordinates.
\epf

Now we immediately see Schlesinger's equations \eqref{eqn: 3x3 schles} become:
\begin{equation}  \label{eqn: xy schles}
\frac{dx}{dt}= \frac{f(x,y)}{t-1},\quad
\frac{dy}{dt}= \frac{-f(x,y)}{t}
\end{equation}
where $f(x,y)=\tr(B_1[B_2,B_3])$.
To write this in terms of $x,y$ we note $f=p-q$ so that
$$f^2 = (p+q)^2 -4pq = (ax+by+k)^2 +4xy(x+y-c)$$
and therefore $f$ is the square root of a cubic polynomial.
The following is then immediate:
\begin{lem}\label{lem transn}
By translating $x,y$ by constants, any cubic polynomial of the form
$$4x^2y+4y^2x + \text{quadratic terms}$$
maybe put in the form
$$4x^2y+4y^2x + Axy+Bx+Cy+D$$
for constants $A,B,C,D$.
\end{lem}

Now in \cite{Hit-GASE} Hitchin has carried out the analogous procedure for the
rank two (trace-free) 
Schlesinger equations \eqref{eqn: 2x2 schles} (which we know are
equivalent to PVI) using local coordinates $x:=\tr(A_1A_3)$ and
$y:=\tr(A_2A_3)$ on the fibre of \eqref{eqn: rk 2+ bundle}, 
and he found that they become equation 
\eqref{eqn: xy schles} but with $f$ replaced by the function
$f_{\text{Hitchin}}$ which satisfies 
\begin{equation}\label{eqn: hitchins f^2}
f_{\text{Hitchin}}^2=-2\det\left(
\begin{matrix}
\epsilon_1 & \epsilon-x-y & x \\
\epsilon-x-y & \epsilon_2 & y \\
x & y & \epsilon_3 
\end{matrix}
\right)
\end{equation}
where $\epsilon_i:=\tr(A_i^2)$ and 
$2\epsilon:=\epsilon_4-\epsilon_1-\epsilon_2-\epsilon_3$ are
constants. 
One easily sees this is again a cubic polynomial of the form appearing in
Lemma \ref{lem transn}; This shows that the system of equations 
\eqref{eqn: xy schles} is equivalent to PVI, and that the four constants
$A,B,C,D$ parameterising the cubic polynomials correspond to the four
parameters $(\al,\be,\ga,\delta)$ of PVI. 
In turn we deduce the desired result 
that the rank three Schlesinger equations \eqref{eqn: 3x3 schles} are also
equivalent to PVI.

Finally let us relate the parameters of the $3\times 3$
 systems to those of the
corresponding $2\times 2$ systems and in turn to the parameters 
$\al,\be,\ga,\delta$ of 
PVI. Define $\mu_1,\mu_2,\mu_3\in \IC$ 
such that the eigenvalues of $B_4\in O_4$
are $\{-\mu_1,-\mu_2,-\mu_3\}$. 
Then the parameters of the $3\times 3$ system are 
$\{\lambda_i\},\{\mu_i\}$ for $i=1,2,3$ 
and (by taking the trace of $\sum B_i=0$) we see
they are constrained by $\sum \lambda_i=\sum\mu_i$.
On the other hand if we define $\theta_i$ to be the difference
between the eigenvalues (in some order) of $A_i$ for $i=1,2,3,4$ then the 
parameters for the $2\times 2$ system are $\{\theta_i\}$ and
the parameters used by Hitchin are $\epsilon_i=\tr(A^2_i)=\theta_i^2/2$.
From \cite{JM81} the relation between $(\al,\be,\ga,\delta)$ 
and  $\{\theta_i\}$ is
\begin{equation}\label{thal}
\al=(\theta_4-1)^2/2, \ \be=-\theta_1^2/2, \ 
\ga=\theta^2_2/2, \ \delta=(1-\theta_3^2)/2.
\end{equation}

To obtain $\{\theta_i\}$ from $\{\lambda_i\},\{\mu_i\}$ we should go via
the parameters $A,B,C,D$ (using $f_{\text{Hitchin}}$ and $f$ respectively).
This appears to be difficult since in either case $A,B,C,D$
are complicated degree six 
polynomials in $\{\theta_i\}$ or $\{\lambda_i,\mu_i\}$
respectively. 
However in fact there are {\em linear} maps from 
$\{\lambda_i,\mu_i\}$ to $\{\theta_i\}$ leading to corresponding parameters
$A,B,C,D$:
\begin{lem} \label{lem params}
If we define
\begin{equation}\label{thla}
\theta_i:= \lambda_i-\mu_1\  (i=1,2,3),\quad\theta_4:=\mu_2-\mu_3
\end{equation}
and then define $ f, f_{\text{Hitchin}}$ as above in terms of 
 $\{\lambda_i, \mu_i\}$, $\{\theta_i\}$ respectively then
$$f^2(x,y)=f^2_{\text{Hitchin}}(x-\theta_1\theta_3/2,y-\theta_2\theta_3/2).$$
Moreover 
the same result holds if $\mu_1,\mu_2,\mu_3$ are permuted arbitrarily.
\end{lem} 
\pf
This may be proved by direct calculation.
\epf

Thus the parameters of PVI corresponding to $\{\lambda_i,\mu_i\}$ are 
given immediately by combining \eqref{thal} and \eqref{thla}. (The other $5$ 
permutations of $\{\mu_1,\mu_2,\mu_3\}$ give equivalent parameters.)

\begin{rmk}
If we choose $B_4$ to be diagonal then 
each of the six off-diagonal entries of the matrix
$$z(z-1)(z-t)B(z)$$ 
of polynomials in $z$, 
is  linear, 
and so has a unique zero on the complex plane.
One may then prove (analogously to the $2\times 2$ 
case) 
that the positions $\eta_{ij}$ of each of these zeros (as functions of $t$)
are solutions of 
PVI. (The parameters of PVI for each of the six solutions $\eta_{ij}$
correspond to one of the six permutations of 
 $\{\mu_1,\mu_2,\mu_3\}$ in Lemma \ref{lem params}).
\end{rmk}

\end{section}
\begin{section}{Braid group actions}

Now we consider the spaces of monodromy data corresponding to the above rank
three systems and describe the natural braid group actions on them.

Let $\cC_i:= \exp(2\pi\sqrt{-1}\cO_i)\subset GL(V)$ be the conjugacy class
associated to the fixed adjoint orbit $\cO_i$.
For $i=1,2,3$ let $t_i:=\exp(2\pi\sqrt{-1}\lambda_i)$ so that $\cC_i$ is four
dimensional and contains pseudo-reflections of the form
$$r_i=1+e_i\otimes\al_i\qquad\text{where}\quad e_i\in V, \al_i\in V^* 
\text{ and }
1+\al_i(e_i)=t_i\ne 0.$$
We suppose that $\cO_4$ is sufficiently generic that $\cC_4$ is a generic
conjugacy class (i.e. that the difference between any two eigenvalues of any
element of $\cO_4$ is not an integer).

Now the monodromy representation of the Fuchsian system 
\eqref{eqn: 3x3 linear} is an element
$$\rho\in \Hom(\pi_1(\IC\setminus\{a_i\},p),G)$$
where $p$ is a base point. If $B_i\in \cO_i$ then $\rho$ maps a simple
positive loop around $a_i$ into $\cC_i$.
Clearly $\rho$ is determined by its values
$$(r_1,r_2,r_3):=(\rho(\ga_1),\rho(\ga_2),\rho(\ga_3))$$
on any set of loops $\{\ga_i\}$ generating $\pi_1(\IC\setminus\{a_i\})$.
By sliding the points $a_1,a_2,a_3$ around the complex plane one obtains an
action of the three-string braid group $B_3$ on the set of such triples:
Two generators of this action are:
$$\be_1(r_1,r_2,r_3) = (r_2,r_2^{-1}r_1r_2,r_3),\qquad 
\be_2(r_1,r_2,r_3) = (r_1,r_3,r_3^{-1}r_2r_3).$$
(One may think of this action as choosing different generators of 
$\pi_1(\IC\setminus\{a_i\})$ as the loops are dragged around when the points
$a_i$ are permuted in the plane.)

Now we wish to descend this action to the space 
of conjugacy classes of fundamental
group representations. 
Let $\{e_i\}$ be the standard basis of $V$.
If each $r_i$ is  pseudo-reflection and $\rho$ is sufficiently generic
then the triple $(r_1,r_2,r_3)$ is conjugate to a triple with 
$r_i=1+e_i\otimes\al_i$ for some $\al_i\in V^*$. 
If we let $u_{ij}=\al_i(e_j)$ then the corresponding
 matrix $U$ (whose rows represent the linear forms $\al_i$) 
is determined by the conjugacy class of
$\rho$ up to conjugation by a diagonal matrix. Thus we have the following five
functions on the space $\Hom_{\text{p-r}}(\pi_1(\IC\setminus\{a_i\}),G)/G$
of conjugacy classes of representations having pseudo-reflection monodromy
around each $a_i$:
$$w=u_{12}u_{21},\  x=u_{13}u_{31},\  y=u_{23}u_{32},\  
p=u_{12}u_{23}u_{31},\  q=u_{32}u_{21}u_{13},$$
which are the ``multiplicative analogues'' of the functions with the same
labels in the previous section (and similarly we see they generate
the ring of $G$ invariant functions on
$\Hom_{\text{p-r}}(\pi_1(\IC\setminus\{a_i\}),G)$). 
The functions $w,x,y,p,q$ may be expressed directly in
terms of the $r_i$ using the formulae
$$\tr(r_ir_j)= 1+t_i+t_j+u_{ij}u_{ji},$$
$$ 
  \tr(r_1r_2r_3)=t_1+t_2+t_3+w+x+y+p, \quad 
  \tr(r_3r_2r_1)=t_1+t_2+t_3+w+x+y+q.$$ 
It is now straightforward to calculate the induced $B_3$ action on the matrix
  $U$ and in turn on the quintuple of functions $w,x,y,p,q$:
If we assume that each $r_i$ is an order two complex reflection
  (i.e. $t_i=-1$) then the
  formula simplifies to
$$\be_1(w,x,y,p,q)=(w,y+p+q+wx,x,-q-wx,-p-wx)$$
$$\be_2(w,x,y,p,q)=(x+p+q+wy,w,y,-q-wy,-p-wy)$$
(and is not much more complicated in general, but this case is sufficient for
the example we will consider here).

Now if we consider the map 
$$ \IP^1\setminus\{0,1,\infty\}\to\IC^3\setminus\text{diagonals},\quad
t\mapsto(a_1,a_2,a_3)=(0,t,1)$$
then loops around $0,1$ generating the fundamental group of
the three-punctured sphere map to the generators 
$\be_1^2,\be_2^2$ of the pure braid group 
$P_3=\pi_1(\IC^3\setminus\text{diags})$ (i.e. to the squares of the chosen 
generators of $B_3$).
From the picture sketched in the introduction it follows that 
the branching of solutions to
the sixth Painlev\'e equation PVI are given by this action of $P_3$. 
Therefore finite orbits of the $P_3$ action correspond to finite branching
solutions to PVI.
In particular if $(r_1,r_2,r_3)$ are a triple of complex reflections (finite
order pseudo-reflections) that generate a finite group then we can be sure to
obtain a finite $P_3$ orbit. 
Thus each triple of generators of each three-dimensional complex reflection
group gives a finite branching solution to PVI.

\begin{rmk}
One may determine the parameters $\{\lambda_i,\mu_i\}$ from the triple 
$(r_1,r_2,r_3)$
(and so by Lemma \ref{lem params} the parameters of PVI) using the fact
that $r_i$ has eigenvalues $\{1,1,e^{2\pi \sqrt{-1}\lambda_i}\}$
and that the product $r_1r_2r_3$ has eigenvalues 
$\{e^{2\pi i\mu_1},e^{2\pi i\mu_2},e^{2\pi i\mu_3}\}$.
In particular by 5.4 of \cite{Shep-Todd} 
if $(r_1,r_2,r_3)$ are one of the standard triples of generators
of a finite complex reflection group $G$, 
then $\mu_i$ are related to the exponents 
$m_1\leqslant m_2\leqslant m_3$ of $G$  
as follows:
$$ \mu_i= m_i/h,\quad (i=1,2,3)\qquad h:=m_3+1.$$ 
This enables us to compile Table \ref{solution table} of 
parameters for solutions from standard generating triples, 
where a suitable permutation of $\{\mu_i\}$ has been used in each case 
and we have taken each $\theta_i$ to be positive (since negating any
$\theta_i$ leads to equivalent PVI parameters).
\end{rmk}

\begin{table}
\begin{center}
\begin{tabular}[h]{|c|c|c|}
\hline
Group & Degrees $m_i+1$ & $(\theta_1,\theta_2,\theta_3,\theta_4)$ \\ \hline
\hline
$G(m,m,3)$ & $3,m,2m$  & $(m-2,m-2,m-2,m)/2m$ \\ \hline
$G(m,1,3)$ & $m,2m,3m$  & $(m-2,m-2,2m-4,4m)/6m$ \\ \hline
Icosahedral & $2,6,10$  & $(0,0,0,4/5)$ \\ \hline
$G_{336}$ & $4,6,14$  & $(2,2,2,4)/7$ \\ \hline
$G_{648}$ & $6,9,12$  & $(0,0,0,1/2)$ \\ \hline
$G_{1296}$ & $6,12,18$  & $(4,7,7,12)/18$ \\ \hline
$G_{2160}$ & $6,12,30$  & $(5,5,5,9)/15$ \\ \hline
\end{tabular}
\vspace{0.2cm}
\caption{Parameters for solutions from standard generating triples.}
\label{solution table}
\end{center}\end{table}

\end{section}

\begin{section}{An example}

As an example let us consider the smallest exceptional 
three dimensional non-real
complex reflection group $K\subset G$ of order $336$. 
(The associated collineation group in $PGL_3(\IC)$ is Klein's simple group of
order $168$.)
To classify the finite branching solutions of PVI associated to $K$ we must
classify up to conjugacy the braid group orbits of generating triples of
reflections in $K$. (Each such triple will give a solution but conjugate
triples, and those in the same braid group orbit, will give the same solution.)
Let $(r_1,r_2,r_3)$ denote the standard triple of generators of $K$ (given as
explicit $3\times 3$ matrices on p.295 of \cite{Shep-Todd}).
Now $K$ contains precisely $21$ complex
reflections all of which have order $2$ and are all conjugate in $K$.
Thus
it is sufficient to consider braid orbits of elements of the form
$$ (r_1,a,b)$$
for reflections $a,b$.
There are $441=21^2$ such triples but it turns out (using Maple) 
that these constitute just
$45$ distinct conjugacy classes of triples 
(i.e. the quintuple of functions $(w,x,y,p,q)$ takes 
only $45$ distinct values on these $441$ triples). 
Thus there are $45$ conjugacy classes of triples of reflections in $K$, since
every reflection is conjugate to $r_1$.
Then it is quite manageable to calculate the braid group orbits on these
conjugacy classes of triples: one finds the $45$ classes are partitioned into
orbits of size
$$1,1,3,3,4,4,6,7,7,9.$$ 
Then we find that, except for the orbits of size $7$, all the corresponding
triples generate proper subgroups of $K$. On the other hand the two orbits of
size seven come from the triples
$$(r_1,r_2,r_3),\qquad\text{and}\qquad(r_1,r_3,r_2)$$
of generators of $K$.
Let us focus on the orbit of the first triple $(r_1,r_2,r_3)$ (the other 
triple gives a solution to PVI with equivalent parameters).
One finds this orbit does
not break up into smaller orbits up when we restrict to the pure braid group;
the $P_3$ orbit still has size seven and so
by the remarks of the preceding sections this implies the existence of a
solution to PVI with seven branches. 
From Table \ref{solution table} this solution has $\theta$-parameters
$(2,2,2,4)/7$.
By examining the permutation of the branches at each of the three
branch-points (i.e. the cycle decomposition of the action of the 
two generators of $P_3$  and of their product) we find
this solution is single-valued on a genus zero covering (at each branch-point
one finds a 3-cycle and two 2-cycles). 
To the authors knowledge such a solution
does not appear in the existing literature.

\end{section}

\begin{section}{Conclusion/Outlook}

In summary 
we have shown how the general PVI equation governs isomonodromic deformations
of rank three systems having four singularities on the sphere with
rank one residue at three of the singularities.
The corresponding space of monodromy data consists of representations of the
fundamental group of the four-punctured sphere such that three of the local
monodromies are pseudo-reflections.
It follows that the branching of the solutions to PVI (i.e. the ``nonlinear
monodromy of PVI'') is governed by the action of the pure three-string braid
group on the set of conjugacy classes of such representations.
Thus for each finite subgroup of $GL_3(\IC)$ generated by three
pseudo-reflections (i.e. each triply-generated three-dimensional complex
reflection group) 
we obtain a finite braid group orbit and thus a finite branching solution to
PVI.
Finally we have started to describe some of the new solutions that arise in
this way.

Some remaining questions are:

\noindent
1) Are all these solutions algebraic? (It is hard to believe they are not,
and a proof should be possible as in \cite{DubMaz00} using/adapting 
Jimbo's study \cite{Jimbo82} of the asymptotics of PVI.)
\footnote{
By finding the corresponding $2\times 2$ triples of monodromy data 
we have recently found that
in fact Jimbo's work may be used directly to prove 
that the solutions found in Section 4 are indeed  
algebraic.
Details will appear elsewhere.}

\noindent
2) What are the explicit formulae for the solutions?
(This could be amenable to brute force methods on a computer, at least for
the smaller braid group orbits, once the asymptotics in 1) are understood.
However, the fact that the icosahedral solution with 18 branches took
9 pages to write down implicitly in the preprint version of \cite{DubMaz00},
leads us to question how valuable such explicit formulae are.)

\noindent
3) Is there a  geometrical or physical interpretation of these solutions? (This
   appears to be a deeper question: for example i) Dubrovin (see \cite{DubPT})
   has shown how
   solutions to $\text{PVI}_{\mu}$ may be used to construct three-dimensional 
semisimple Frobenius manifolds (i.e. certain approximations to $2\text{D}$
   topological quantum field theories), and ii) Hitchin  
   \cite{Hit-Poncelet, Hit-Ed} and Doran
   \cite{Chuck1} have constructed some of the algebraic solutions purely
   geometrically relating them to previously solved, often classical, 
   algebro-geometric problems.)

\end{section}

\renewcommand{\baselinestretch}{1}              %
\normalsize
\bibliographystyle{amsplain}    \label{biby}
\bibliography{../thesis/syr}    

\providecommand{\bysame}{\leavevmode\hbox to3em{\hrulefill}\thinspace}
\providecommand{\MR}{\relax\ifhmode\unskip\space\fi MR }
% \MRhref is called by the amsart/book/proc definition of \MR.
\providecommand{\MRhref}[2]{%
  \href{http://www.ams.org/mathscinet-getitem?mr=#1}{#2}
}
\providecommand{\href}[2]{#2}
\begin{thebibliography}{10}

\bibitem{smid}
P.P. Boalch, \emph{{S}ymplectic manifolds and isomonodromic deformations}, Adv.
  in Math. \textbf{163} (2001), 137--205,
  {(http://www.idealibrary.com/links/toc/aima/163/2/0)}.

\bibitem{bafi}
\bysame, \emph{G-bundles, isomonodromy and quantum {W}eyl groups}, Int. Math.
  Res. Not. (2002), no.~22, 1129--1166, math.DG/0108152.

\bibitem{BMM99}
M.~Brou{\'e}, G.~Malle, and J.~Michel, \emph{Towards spetses. {I}}, Transform.
  Groups \textbf{4} (1999), no.~2-3, 157--218, Dedicated to the memory of
  Claude Chevalley. \MR{2001b:20082}

\bibitem{DKP}
C.~De~Concini, V.~G. Kac, and C.~Procesi, \emph{Quantum coadjoint action}, J.
  Amer. Math. Soc. \textbf{5} (1992), no.~1, 151--189.

\bibitem{Chuck1}
C.~F. Doran, \emph{Algebraic and geometric isomonodromic deformations}, J.
  Differential Geom. \textbf{59} (2001), no.~1, 33--85. \MR{1 909 248}

\bibitem{DubPT}
B.~Dubrovin, \emph{Painlev\'e transcendents in two-dimensional topological
  field theory}, The Painlev\'e property, Springer, New York, 1999,
  pp.~287--412, (math/9803107).

\bibitem{DubMaz00}
B.~Dubrovin and M.~Mazzocco, \emph{Monodromy of certain {P}ainlev\'e-{V}{I}
  transcendents and reflection groups}, Invent. Math. \textbf{141} (2000),
  no.~1, 55--147. \MR{2001j:34114}

\bibitem{Hit95long}
N.~J. Hitchin, \emph{Frobenius manifolds}, Gauge Theory and Symplectic Geometry
  (J.~Hurtubise and F.~Lalonde, eds.), NATO ASI Series C: Maths \& Phys., vol.
  488, Kluwer, 1995.

\bibitem{Hit-Poncelet}
\bysame, \emph{Poncelet polygons and the {P}ainlev\'e equations}, Geometry and
  analysis (Bombay, 1992), Tata Inst. Fund. Res., Bombay, 1995, pp.~151--185.
  \MR{97d:32042}

\bibitem{Hit-GASE}
\bysame, \emph{Geometrical aspects of {S}chlesinger's equation}, J. Geom. Phys.
  \textbf{23} (1997), no.~3-4, 287--300. \MR{99a:32023}

\bibitem{Hit-Ed}
\bysame, \emph{Quartic curves and icosahedra}, talk at Edinburgh, September
  1998.

\bibitem{Jimbo82}
M.~Jimbo, \emph{Monodromy problem and the boundary condition for some
  {P}ainlev\'e equations}, Publ. Res. Inst. Math. Sci. \textbf{18} (1982),
  no.~3, 1137--1161. \MR{85c:58050}

\bibitem{JM81}
M.~Jimbo and T.~Miwa, \emph{Monodromy preserving deformations of linear
  differential equations with rational coefficients {II}}, Physica 2D (1981),
  407--448.

\bibitem{KPS01}
L.~Katzarkov, T.~Pantev, and C.~Simpson, \emph{Density of monodromy actions on
  non-abelian cohomology}, math.AG/0101223.

\bibitem{Shep-Todd}
G.~C. Shephard and J.~A. Todd, \emph{Finite unitary reflection groups},
  Canadian J. Math. \textbf{6} (1954), 274--304. \MR{15,600b}

\bibitem{tot}
B.~Totaro, \emph{Towards a {S}chubert calculus for complex reflection groups},
  Math. Proc. Camb. Phil. Soc., to appear,
  www.dpmms.cam.ac.uk/~bt219/hall.dvi.gz.

\end{thebibliography}
\ 

\ 

\ 

\ 

[Journal ref: Ann. Inst. Fourier, 53, 4 (2003) 1009-1022. ]

\end{document}